\documentclass[12pt,english]{amsart}
\usepackage[cp1251]{inputenc}
\usepackage[english]{babel}
\usepackage{amsmath,amsthm,amssymb}
\newtheorem{theorem}{Theorem}[section]
\newtheorem{lemma}[theorem]{Lemma}
\newtheorem{corollary}[theorem]{Corollary}

\numberwithin{equation}{section}

\begin{document}

\title[Total variation distance estimates]
{Total variation distance estimates via $L^2$-norm
for polynomials in log-concave random vectors}

\author{Egor D. Kosov}

\maketitle

\begin{abstract}
The paper provides an estimate
of the total variation distance between
distributions of polynomials
defined on a space equipped with a logarithmically
concave measure in terms of the $L^2$-distance
between these polynomials.
\end{abstract}

\noindent
Keywords:
logarithmically concave measure, Gaussian measure,
total variation distance,
distribution of a polynomial

\noindent
AMS Subject Classification: 60E05, 60E15, 28C20, 60F99

\section{Introduction}
Davydov and Martynova \cite{DavM} formulated the following
interesting property of polynomials on a space with a Gaussian measure.

{\bf Theorem A.}
{\it
Let $d\in\mathbb{N}$ and let $g$ be a non-constant polynomial of degree $d$ on $\mathbb{R}^n$.
Then there is a constant $C(d, g)$
depending only on $d$ and $g$ such that for every polynomial $f$ of degree $d$
one has
$$
\|\gamma\circ f^{-1} - \gamma\circ g^{-1}\|_{\rm TV}\le
C(d, g) \|f-g\|_{L^2(\gamma)}^{1/d},
$$
where $\gamma$ is the standard Gaussian measure on $\mathbb{R}^n$ and
$\gamma\circ f^{-1}$ and $\gamma\circ g^{-1}$
are the distributions of random variables $f$ and $g$, respectively.
}

Note that in \cite{DavM} the assertion was formulated in terms of multiple stochastic
integrals of order $d$, but the claim above is equivalent to the original one.
The cited paper contains no technical details of the proof,
which, however, can be found in Martynova's PhD thesis.
Nevertheless, since these details are still unpublished and hardly accessible
(Martynova's PhD thesis can be only found in some libraries in Saint Petersburg and Moscow),
there have been several attempts to give a full proof to the above result.
Firstly, Nourdin and Poly \cite{NP} obtained the following theorem.

{\bf Theorem B.}
{\it
Let $d\in\mathbb{N}$, $a>0$, $b> 0$.
Then there exists a number $C(d, a, b)>0$
such that for every pair of polynomials $f, g$ of degree~$d$ on $\mathbb{R}^n$
one has
$$
\|\gamma\circ f^{-1} - \gamma\circ g^{-1}\|_{\rm TV}
\le C(d, a, b)\|f-g\|_{L^2(\gamma)}^{1/(2d)},
$$
provided that the variance of $g$ is in $[a,b]$.
}

The above theorem
clarifies some dependence of $C(d, g)$ on~$g$:
it depends only on the bounds for the variance.
However,
the power of the $L^2$-norm in the theorem is
worse than in Theorem A.
Next, in \cite{BKZ} an intermediate result between Theorem A and Theorem B
was obtained. The constant there was worse than in the Nourdin--Poly estimate,
but the dependence on the $L^2$-norm
differed from the one in \cite{DavM}  by only a logarithmic factor.
Finally, in \cite{Zel} the following theorem was proved.

{\bf Theorem C.}
{\it
Let $d\in\mathbb{N}$.
There is a constant $c(d)$ depending only on $d$ such that,
for every pair of polynomials $f, g$ of degree $d\ge2$ on $\mathbb{R}^n$,
one has
$$
\|\gamma\circ f^{-1} - \gamma\circ g^{-1}\|_{\rm TV}\le
c(d)\bigl(\|\nabla g\|_*^{-1/(d-1)} + 1\bigr)\|f-g\|_2^{1/d},
$$
where
$$
\|\nabla g\|_{*}^2:= \sup\limits_{|e|=1}\int |\partial_e g|^2\, d\gamma.
$$
}

Note that while Theorem C coincides with the Davydov--Martynove estimate,
the constant there is still worse than in the Nourdin--Poly estimate.

This paper generalizes the Davydov--Martynove bound
to the case of an arbitrary log-concave measure in place of a Gaussian measure.
Recall that a probability Borel measure $\mu$ on $\mathbb{R}^n$ is called
logarithmically concave (log-concave or convex) if
$$
    \mu(t A + (1-t)B) \ge \mu(A)^{t}\mu(B)^{1-t} \quad \forall\, t\in [0,1]
$$
for all Borel sets $A,B\subset \mathbb{R}^n$ (see \cite{Bor} and discussion in
\cite[Section 3.10(vi)]{MeasTh} and in \cite[Section 4.3]{DiffMeas}).
This is equivalent to the fact that the measure~$\mu$
has a density of the form $e^{-V}$ with respect to Lebesgue measure
on some affine subspace~$L$,
where $V\colon\, L\to (-\infty, +\infty]$ is a convex function.
We also recall that the total variation norm of a (signed) measure $\nu$ on $\mathbb{R}^n$
is defined by the equality
$$
\|\nu\|_{\rm TV} := \sup\biggl\{\int \varphi\, d\nu, \ \varphi\in C_0^\infty(\mathbb{R}^n),\ \|\varphi\|_\infty\le1\biggr\},
$$
where $\|\varphi\|_\infty:=\sup|\varphi(x)|$.
The distribution $\mu\circ F^{-1}$
of a measurable function $F$ on a measurable space equipped with
a measure $\mu$ is a measure on the real line such that
$\mu\circ F^{-1}(A):=\mu(F\in A)$ for all Borel sets $A$.

The main result of the present paper asserts that,
for an all $n,d\in\mathbb{N}$, $d\ge2$,
there is a constant $C(d)$ such that,
for every log-concave measure $\mu$, and
every pair of polynomials $f,g$ of degree $d$
on $\mathbb{R}^n$,
one has
$$
\sigma_g^{1/d}\|\mu\circ f^{-1} - \mu\circ g^{-1}\|_{\rm TV}\le
C(d)\|f-g\|_2^{1/d},
$$
where
$\sigma_g^2:=\mathbb{D}g$ is the variance of $g$.
We note that even in the case of a Gaussian measure the obtained result
improves the dependence of the constant in comparison to Theorem~C.
We also note that due to independence of the constant in the inequality
of the dimension the same estimate remains valid in the
infinite dimensional case.
The proof of the announced inequality
develops some ideas of \cite{Kos}, \cite{NP}, and \cite{Zel}.

\section{Preliminaries}

This section contains necessary definitions, notation,
and several known results which are used further.
We mainly consider the finite-dimensional space
$\mathbb{R}^n$ equipped with
the Borel $\sigma$-field and with the
standard Euclidian inner product $(x, y)$, $x,y\in \mathbb{R}^n$.
Let $|\cdot|$ be the standard norm
$|x|:= \sqrt{(x,x)}$, $x\in \mathbb{R}^n$.
Let $C_0^\infty(\mathbb{R}^n)$ denote the space
of all infinitely
smooth functions with compact support.

A log-concave measure $\mu$ on $\mathbb{R}^n$ is called isotropic
if it is absolutely continuous with respect to Lebesgue measure and
$$
\int_{\mathbb{R}^n}(x,\theta)\, \mu(dx)=0,\quad
\int_{\mathbb{R}^n}(x,\theta)^2\, \mu(dx)=
|\theta|^2\quad \forall\ \theta\in\mathbb{R}^n.
$$

The Skorohod derivative $D_e\mu$ of a Borel
measure $\mu$ along a vector
$e\in\mathbb{R}^n$
is a bounded signed Borel measure on $\mathbb{R}^n$ such that
$$
\int_{\mathbb{R}^n} \partial_e \varphi\, d\mu=
-\int_{\mathbb{R}^n} \varphi\, d(D_e\mu)
$$
for every $\varphi\in C_0^\infty(\mathbb{R}^n)$ (see \cite{DiffMeas}).
It was proved by Krugova \cite{Krug}
(see also \cite[Section~4.3]{DiffMeas}) that
every log-concave measure $\mu$ with density $\rho$
is Skorohod differentiable along every vector $e\in\mathbb{R}^n$
and for every unit
vector $e$ one has
$$
    \|D_e\mu\|_{\rm TV}=2\int_{\langle e\rangle^\bot} \max\limits_t\rho(x+te)dx,
$$
where $\langle e\rangle^\bot$ is the orthogonal complement of $e$.

If the measure $\mu$ is fixed,
for a $\mu$-measurable function $f$ we set
$$
\|f\|_r:=\Bigl(\int_{\mathbb{R}^n} |f|^r\, d\mu\Bigr)^{1/r}\ \text{for} \ r>0, \quad
\|f\|_0:=\exp\Bigl(\int_{\mathbb{R}^n} \ln|f|\, d\mu\Bigr)
=\lim\limits_{r\to 0} \|f\|_r.
$$
An important feature of the $0-$``norm'' is its multiplicative property, i.e.
$\|f\cdot g\|_0=\|f\|_0\cdot\|g\|_0$.
We also denote the expectation
and the variance of the random variable $f$
by the symbols $\mathbb{E}f$ and~$\sigma^2_f$ respectively, i.e.
$$
\mathbb{E}f:=\int_{\mathbb{R}^n} f\, d\mu\quad
\sigma^2_f := \mathbb{D}f = \int_{\mathbb{R}^n} (f - \mathbb{E}f)^2\, d\mu.
$$
Throughout the paper the symbols $c, C, c_1, C_1, \ldots$
denote positive universal constants,
the symbols $c(d), C(d), c_1(d), C_1(d),\ldots$
denote positive constants, that depend only on one parameter~$d$,
and $c(d,n), C(d,n), c_1(d,n), C_1(d,n),\ldots$
denote positive constants, that depend only on two parameters $d$ and $n$.
The values of these constants are not necessarily the same in different
appearances.
Throughout the paper we
omit the indication of $\mathbb{R}^n$
in all integrations.

\vskip .1in

We now formulate some key known
results which will be applied
in the proofs.

The first result is the so-called Carbery--Wright inequality
for polynomials on a space with a log-concave measure.
\begin{theorem}[\cite{CarWr}, \cite{NSV}]\label{t1.1}
There is an absolute constant $c$ such that for every log-concave measure
$\mu$ on $\mathbb{R}^n$ and every polynomial
$f$ of degree $d$ the following inequality holds true:
$$
\mu(|f|\le t)\left(\int|f|d\mu\right)^{1/d}\le cd\, t^{1/d} .
$$
\end{theorem}

The next result shows that for a log-concave measure
all the $L^p$-``norms'' on the space of
polynomials of a fixed degree are equivalent.
Those ``norms''
estimate each other with constants depending only
on the degree of polynomials.
\begin{theorem}[\cite{BobkPoly}, \cite{Bobk}]\label{t1.2}
Let $\mu$ be a log-concave measure on $\mathbb{R}^n$, $q\ge1$.
Then there is an absolute constant $c$ such that for every
polynomial $f$ of degree $d$ the following inequalities hold:
$$
\|f\|_q\le (cqd)^d\|f\|_0, \quad \|f\|_q\le (cq)^d\|f\|_1.
$$
\end{theorem}

We also need the following results
on the structure of the density of a log-concave measure.
The next theorem can be found in \cite[Proposition 4.1]{BobkIsop}.
\begin{theorem}[\cite{BobkIsop}]\label{t1.3}
Let $\mu$ be a log-concave measure on $\mathbb{R}$ with density $\rho$.
Then
$$
\|\rho\|^2_\infty\int\biggl(t - \int\tau\mu(d\tau)\biggr)^2\mu(dt)\ge 12^{-1}.
$$
\end{theorem}

\begin{theorem}[\cite{Klartag}, \cite{Ball}]\label{t1.8}
For every $n\in \mathbb{N}$, there is a constant $C(n)$ depending only on $n$ such that
for every isotropic log-concave measure $\mu$ on $\mathbb{R}^n$
with density $\rho$
one has
$$
(\max\rho)^{1/n}\le C(n).
$$
\end{theorem}

There is a conjecture
(the hyperplane conjecture) that the constant above
can be chosen independent of $n$, but
the best known constant so far is  $C_n\sim n^{1/4}$,
which is due to Klartag \cite{Klartag}.

The following result is Corollary 2.4 in \cite{Klartag07}.
\begin{theorem}[\cite{Klartag07}]\label{t1.4}
For every $n\in\mathbb{N}$,
there are universal constants $C, c>0$
such that, for every
isotropic log-concave measure $\mu$ on $\mathbb{R}^n$
with density $\rho$,
the following inequality holds:
$$
\rho(x)\le \rho(0) e^{Cn - c|x|}.
$$
\end{theorem}

The next property is a combination of Corollary 5.3
and Lemma 5.4 in \cite{Klar-CLT}.
\begin{theorem}[\cite{Klar-CLT}]\label{t1.5}
Let $n\in\mathbb{N}$, $n\ge2$ and let $\mu$ be an
isotropic log-concave measure on $\mathbb{R}^n$
with density $\rho$. Let
$K=\{x\in\mathbb{R}^n: \rho(x)\ge e^{-20n}\rho(0)\}$.
Then
$$
B_\frac{1}{10}\subset K.
$$
\end{theorem}

The following theorem states
the Poincar$\acute{e}$ inequality for log-concave measures.
\begin{theorem}[\cite{BobkIsop, KLS}]\label{t1.6}
There is an absolute constant $M$ such that for every log-concave measure
$\mu$ on $\mathbb{R}^n$ and every locally Lipschitz function $f$
one has
$$
\int\left(f-\int fd\mu\right)^2d\mu\le M\int|x-x_0|^2d\mu\int|\nabla f|^2d\mu,
$$
where $x_0=\int x d\mu$.
\end{theorem}

The following so-called localization lemma from \cite{FrGue}
(see also \cite{KLS} and \cite{LS}) plays a crucial role in our proof.

\begin{theorem}[Localization lemma with $p$ constraints, see \cite{FrGue}]\label{t1.7}
Let $K$ be a compact convex set in $\mathbb{R}^n$,
$F_i\colon K\to \mathbb{R}$, $1\le i\le p$.
Assume that all functions $F_i$ are upper semi-continuous.
Let $P_{F_1,\ldots, F_p}$ be the set of
all log-concave measures with support in $K$
such that
$$
\int F_i\, d\mu\ge0, \ i=1,\ldots,p.
$$
Let $\Phi\colon P(K)\to\mathbb{R}$ be a convex
upper semi-continuous
function, where $P(K)$
is the space of all Borel probability measures supported in $K$
equipped with the weak topology.
Then $\sup\limits_{\mu\in P_{F_1,\ldots, F_p}}\Phi(\mu)$
is attained on log-concave measures $\mu$
such that the smallest affine subspace containing the support of $\mu$
is of dimension at most $p$.
\end{theorem}

\section{Total variation distance estimate}

We start with the following reverse
Poincar$\acute{e}$ inequality
for polynomials on a space with a log-concave measure.
Such estimates are well known for Gaussian measures
due to the equivalence of all Sobolev norms
on the space of all polynomials of a fixed degree (see \cite{GM}).

\begin{theorem}\label{t2.1}
Let $n,d\in\mathbb{N}$. There is
a constant $C(d)$,
which depends only on the degree~$d$, such that,
for each log-concave measure $\mu$
on $\mathbb{R}^n$,
each polynomial $f$ of degree~$d$,
and
each vector~$e$ of unit length, one has
$$
\|\partial_ef\|_2\le C(d)\|D_e\mu\|_{\rm TV}\|f\|_2.
$$
\end{theorem}

\begin{proof}
We firstly consider the one-dimensional case.
By homogeneity we can assume that the polynomial $f$ is of the form
$$
f(t)=\prod_{i=1}^d(t-t_i).
$$
Moreover, without loss of generality we can assume that $\int t\, \mu(dt)=0$.
Using Theorem~\ref{t1.2} we get
\begin{multline*}
\int t^2\, \mu(dt)\int (f'(t))^2\, \mu(dt)\le d\sum_{i=1}^d\int t^2\, \mu(dt)\int\Bigl|\prod_{j\ne i}(t-t_j)\Bigr|^2\, \mu(dt)\\
\le d(2cd)^{2d}\sum_{i=1}^d\int t^2\, \mu(dt)\prod_{j\ne i}\int|t-t_j|^2\, \mu(dt)
=d(2cd)^{2d}\sum_{i=1}^d\int t^2\, \mu(dt)\prod_{j\ne i}\left(\int t^2\, \mu(dt)+|t_j|^2\right)\\
\le d^2(2cd)^{2d}\prod_{i=1}^d\left(\int t^2\, \mu(dt)+|t_i|^2\right)=d^2(2cd)^{2d}\prod_{i=1}^d\int|t-t_i|^2\, \mu(dt)
\le d^2(4c^2d)^{2d}\int f^2\, \mu(dt).
\end{multline*}
Thus,
$$
\sigma(\mu)\|f'\|_2\le (Cd)^d\|f\|_2,
$$
where $\sigma^2(\mu)$ is the variance of $\mu$.
The last bound combined with Theorem \ref{t1.3}
implies
$$
\|f'\|_2\le (Cd)^d\|\rho\|_\infty\|f\|_2,
$$
which is equivalent to the inequality
$$
\|f'\|_1\le (Cd)^d\|\rho\|_\infty\|f\|_1
$$
due to Theorem \ref{t1.2}.

We now proceed to the general case.
Without loss of generality we can assume that $e=e_1$ is the first basis vector.
Set $\tilde{x}:=(x_2,\ldots, x_n)$ and
$$
\tilde{\rho}(x_1, x_2, \ldots, x_n):=
\frac{\rho(x_1, x_2, \ldots, x_n)}{\int\rho(\tau, x_2, \ldots, x_n)\, d\tau}.
$$
Applying the obtained one-dimensional bound and Theorem \ref{t1.2} we get
\begin{multline*}
\|\partial_{e_1}f\|_1^{1/2}\le c(d)\int |\partial_{e_1}f|^{1/2}\rho\, dx\\
= c(d)\int \biggl(\int\rho (\tau, \tilde{x})\, d\tau\biggr)\int |\partial_{e_1}f|^{1/2}
\tilde{\rho}(x_1,\tilde{x})\, dx_1d\tilde{x}\\
\le c(d)\int \biggl(\int\rho (\tau, \tilde{x})\, d\tau\biggr)\biggl(\int |\partial_{e_1}f|
\tilde{\rho}(x_1, \tilde{x})dx_1\biggr)^{1/2}\, d\tilde{x}\\
\le c_1(d) \int \biggl(\int\rho (\tau, \tilde{x})\, d\tau\biggr)
\biggl(\max\limits_t\tilde{\rho}(t, \tilde{x})
\int|f|\tilde{\rho}(x_1, \tilde{x})dx_1\biggr)^{1/2}\, d\tilde{x}\\
\le c_1(d) \biggl(\int\max\limits_t\rho (t, \tilde{x})\, d\tilde{x}\biggr)^{1/2}
\biggl(\int|f|\rho\, dx\biggr)^{1/2}
= c_2(d) \|D_{e_1}\mu\|_{\rm TV}^{1/2}\|f\|_1^{1/2}.
\end{multline*}
The theorem is proved.
\end{proof}

We also need the following technical lemma.
\begin{lemma}\label{lem3.1}
Let $n,d\in\mathbb{N}$, $n\ge2$.
There is a constant $c(d,n)$,
depending only on $d$ and $n$, such that,
for every isotropic
log-concave measure $\mu$ on $\mathbb{R}^n$
with density $\rho$,
every polynomial~$h$ of degree $d$,
and
every unit vector $e\in\mathbb{R}^n$, the following bound holds:
$$
\int_{\langle e\rangle^\bot}\max\limits_{s}\bigl[|h(x+se)|\rho(x+se)\bigr]\, dx\le
c(d,n)\int_{\mathbb{R}^n}|h|\, d\mu.
$$
\end{lemma}
\begin{proof}
By Theorem \ref{t1.4} there is a bound
$\rho(x)\le \rho(0) e^{Cn - c|x|}$ implying
\begin{multline*}
\int_{\langle e\rangle^\bot}\max\limits_{s}\bigl[|h(x+se)|\rho(x+se)\bigr]\, dx
\le
\rho(0)\int_{\langle e\rangle^\bot}\max\limits_{s}\bigl[|h(x+se)| e^{Cn - c|x+se|}\bigr]\, dx
\\
\le
\rho(0)e^{Cn}\int_{\langle e\rangle^\bot}e^{-c_1|x|}\max\limits_{s}\bigl[|h(x+se)| e^{- c_1|s|}\bigr]\, dx.
\end{multline*}
We now note that the function $s\mapsto h(x+se)$ is a polynomial. Thus,
$h(x+se)=a_ds^d+a_{d-1}s^{d-1}+\ldots+a_1s+a_0$, where $a_0,\ldots a_d$ are some functions of variable $x$.
Using this representation we can write
$$
\max\limits_{s}\bigl[|h(x+se)| e^{- c_1|s|}\bigr]\le
\sum_{j=0}^d|a_j|\max\limits_{s}[|s|^je^{- c_1|s|}]
=
\sum_{j=0}^d|a_j|\Bigl(\frac{j}{c_1}\Bigr)^je^{-j}
\le
c_1(d)\sum_{j=0}^d|a_j|.
$$
Since all norms on the space of polynomials of a fixed degree on the real line
are equivalent, there is a constant $c_2(d)$ such that
$$
\sum_{j=0}^d|a_j|\le c_2(d)\int_{\mathbb{R}}|h(x+se)|e^{-c_1|s|}\, ds.
$$
Thus,
\begin{multline*}
\int_{\langle e\rangle^\bot}\max\limits_{s}\bigl[|h(x+se)|\rho(x+se)\bigr]\, dx
\le
\rho(0)e^{Cn}c_3(d)\int_{\langle e\rangle^\bot}e^{-c_1|x|}
\int_{\mathbb{R}}|h(x+se)|e^{-c_1|s|}\, ds\, dx
\\
\le
\rho(0)e^{Cn}c_3(d)\int_{\mathbb{R}^n}|h(y)|e^{-c_1|y|}\, dy.
\end{multline*}
Again, since all norms on the space of polynomials of a fixed degree on $\mathbb{R}^n$
are equivalent, there is a constant $c_4(d,n)$ such that
$$
\int_{\mathbb{R}^n}|h(y)|e^{-c_1|y|}\, dy
\le
c_4(d,n)\int_{B_{\frac{1}{10}}}|h(y)|\, dy.
$$
By Theorem \ref{t1.5}
$$
B_\frac{1}{10}\subset K,
$$
where $K=\{y\in\mathbb{R}^n: \rho(y)\ge e^{-20n}\rho(0)\}$,
which implies
$$
\int_{B_{\frac{1}{10}}}|h(y)|\, dy
\le
\int_{K}|h(y)|\, dy
\le
e^{20n}(\rho(0))^{-1}\int_{K}|h(y)|\rho(y)\, dy
\le
e^{20n}(\rho(0))^{-1}\int_{\mathbb{R}^n}|h(y)|\rho(y)\, dy.
$$
Thus, combining the obtained bounds, we get the announced estimate.
\end{proof}

The following technical lemma
provides an estimate similar to the one stated in the introduction,
but is not dimension free. However, it is the main step in the proof of the general result.
\begin{lemma}\label{lem3.2}
Let $n,d\in\mathbb{N}$, $d\ge2$.
There is a constant $c(d,n)$,
which depends only on $d$ and $n$,
such that, for any isotropic log-concave measure
$\mu$ on $\mathbb{R}^n$,
any polynomials~$f,g$ of degree $d$,
any
function $\varphi\in C^\infty_0(\mathbb{R}^n)$
with $\|\varphi\|_\infty\le1$, and any vector $e$ of unit length,
one has
$$
\|\partial_eg\|_2^{1/d}\int\varphi(f) - \varphi(g)\, d\mu
\le c(d,n)\|f-g\|_1^{1/d}.
$$
\end{lemma}

\begin{proof}
Let $\rho=e^{-V}$ be the density of $\mu$, where $V$ is a convex function.
We firstly consider the case $\rho\in C^{\infty}(\mathbb{R}^n)$, $\rho>0$, and $n\ge2$.
Let
$$
\Phi(t) := \int_{-\infty}^t\varphi(\tau)d\tau.
$$
As in \cite{NP}, \cite{BKZ}, and \cite{Zel}, we use the equality
$$
\partial_eg(\varphi(f) - \varphi(g))=
\partial_e(\Phi(f) - \Phi(g))-(\partial_ef-\partial_eg)\varphi(f).
$$
Thus,
\begin{multline*}
\int\bigl(\varphi(f) - \varphi(g)\bigr)\, d\mu
=
\int\frac{(\partial_eg)^2(\varphi(f) - \varphi(g))}{(\partial_eg)^2+\varepsilon}\, d\mu
+
\varepsilon\int\frac{\varphi(f) - \varphi(g)}{(\partial_eg)^2+\varepsilon}\, d\mu\\
=
\int\frac{\partial_eg\partial_e(\Phi(f) - \Phi(g))}{(\partial_eg)^2+\varepsilon}\, d\mu
-
\int\frac{\partial_eg(\partial_ef-\partial_eg)\varphi(f)}{(\partial_eg)^2+\varepsilon}\, d\mu
+
\varepsilon\int\frac{\varphi(f) - \varphi(g)}{(\partial_eg)^2+\varepsilon}\, d\mu.
\end{multline*}
We now estimate each term separately
starting with the last term. By the Carbery--Wright inequality
(Theorem \ref{t1.1}) one has
(see the proof of Lemma 3.1 in \cite{Kos} or expression (4.4) in \cite{BKZ})
\begin{multline*}
\varepsilon\int\frac{\varphi(f) - \varphi(g)}{(\partial_eg)^2+\varepsilon}\, d\mu\le
2c_1d\Bigl(\int_0^\infty(s+1)^{-2}s^{1/(2d-2)}\, ds\Bigr) \|\partial_eg\|_2^{-1/(d-1)}\varepsilon^{1/(2d-2)}
\\
=
C_1(d)\|\partial_eg\|_2^{-1/(d-1)}\varepsilon^{1/(2d-2)}.
\end{multline*}
For the second term we have
$$
-\int\frac{\partial_eg(\partial_ef-\partial_eg)\varphi(f)}{(\partial_eg)^2+\varepsilon}\, d\mu\le
2^{-1}\varepsilon^{-1/2}\int|\partial_ef-\partial_eg|d\mu
\le
C_2(d)\varepsilon^{-1/2}\|D_e\mu\|_{\rm TV}\|f-g\|_1.
$$
Recall that
$$
\|D_e\mu\|_{\rm TV} = 2\int_{\langle e\rangle^\bot}\max\limits_{s}\rho(x+se)\, dx.
$$
Thus, by Theorem \ref{t1.4},
\begin{equation}\label{est1}
\|D_e\mu\|_{\rm TV}\le 2\rho(0)e^{Cn} \int_{\langle e\rangle^\bot}\max\limits_{s}e^{- c|x+se|}\, dx
\le c_1(n)\rho(0)\le c_2(n),
\end{equation}
where Theorem \ref{t1.8} was applied in the last inequality.
Thus,
$$
-\int\frac{\partial_eg(\partial_ef-\partial_eg)\varphi(f)}{(\partial_eg)^2+\varepsilon}\, d\mu\le
c_3(d,n)\varepsilon^{-1/2}\|f-g\|_1.
$$

Now,
integrating by parts in the first term, we get
\begin{multline*}
\int\frac{\partial_eg\partial_e(\Phi(f) - \Phi(g))}{(\partial_eg)^2+\varepsilon}\, d\mu
\\
=
-\int(\Phi(f) - \Phi(g))\Bigl[\frac{\partial^2_eg}{(\partial_eg)^2+\varepsilon}
-2\frac{(\partial_eg)^2\partial^2_eg}{((\partial_eg)^2+\varepsilon)^2}\Bigr]\, d\mu
-\int(\Phi(f) - \Phi(g))\frac{\partial_eg}{(\partial_eg)^2+\varepsilon}\, d(D_e\mu)
\end{multline*}
Up to factor $3$, the first integral above is estimated by
\begin{multline*}
\int\Bigl|\frac{\partial^2_eg}{(\partial_eg)^2+\varepsilon}\Bigr||f-g|\, d\mu=
\int_{\langle e\rangle^\bot}\int_{\mathbb{R}}
\left|\frac{\partial^2_eg(x+te)}{(\partial_eg(x+te))^2+\varepsilon}\right||f(x+te)-g(x+te)|\rho(x+te)\, dtdx\\
\le d\varepsilon^{-1/2}\int_{\langle e\rangle^\bot}\Bigl(\int_{\mathbb{R}}\Bigl|\frac{1}{\tau^2+1}\Bigr|\, d\tau\Bigr)
\max\limits_{s}\bigl[|f(x+se)-g(x+se)|\rho(x+se)\bigr]\, dx
\\
=\pi d\varepsilon^{-1/2}
\int_{\langle e\rangle^\bot}\max\limits_{s}\bigl[|f(x+se)-g(x+se)|\rho(x+se)\bigr]\, dx
\le
c_4(d,n)\varepsilon^{-1/2}\int|f-g|\, d\mu,
\end{multline*}
where Lemma \ref{lem3.1} was applied in the last inequality.
The second integral is not greater than
\begin{multline*}
\varepsilon^{-1/2}\int|f - g|\, d|D_e\mu|\le
\varepsilon^{-1/2}\sqrt{\|D_e\mu\|_{\rm TV}}\Bigl(\int|f - g|^2\, d|D_e\mu|\Bigr)^{1/2}
\\
\le
\varepsilon^{-1/2}\sqrt{c_2(n)}\Bigl(\int|f - g|^2\, d|D_e\mu|\Bigr)^{1/2}.
\end{multline*}
Since $\rho=e^{-V}$, we have $D_e\mu = -V'_ee^{-V}dx$
and $|D_e\mu|=|V'_e|e^{-V}$. For a point $x\in\langle e\rangle^\bot$
let $T(x)$ be such that $V'_e(x+T(x)e)=0$.
Then $V'_e(x+te)\le0$ for $t<T(x)$ and $V'_e(x+te)\ge0$ for $t>T(x)$
by the convexity of the function $t\mapsto V(x+te)$.
Thus,
\begin{multline*}
\int_{\mathbb{R}}|f(x+te)-g(x+te)|^2|V'_e(x+te)|e^{-V(x+te)}\, dt
\\
=
-\int_{-\infty}^{T(x)}
|f(x+te)-g(x+te)|^2V'_e(x+te)e^{-V(x+te)}\, dt
\\
+
\int_{T(x)}^\infty
|f(x+te)-g(x+te)|^2V'_e(x+te)e^{-V(x+te)}\, dt
\\
=
2|f(x+T(x)e)-g(x+T(x)e)|^2\rho(x+T(x)e)
\\-
2\int_{-\infty}^{T(x)}
(\partial_ef(x+te) - \partial_eg(x+te))(f(x+te)-g(x+te))\rho(x+te)\, dt
\\
+
2\int_{T(x)}^\infty
(\partial_ef(x+te) - \partial_eg(x+te))(f(x+te)-g(x+te))\rho(x+te)\, dt
\\
\le
2\max\limits_{s}\bigl[|f(x+se)-g(x+se)|^2\rho(x+se)\bigr]
\\
+
4\int_{\mathbb{R}}|\partial_ef(x+te) - \partial_eg(x+te)||f(x+te)-g(x+te)|\rho(x+te)\, dt.
\end{multline*}
Therefore, we have
\begin{multline*}
\int|f - g|^2\, d|D_e\mu|
=
\int_{\langle e\rangle^\bot}\int_{\mathbb{R}}
|f(x+te)-g(x+te)|^2|V'_e(x+te)|e^{-V(x+te)}\, dt\, dx
\\
\le
2\int_{\langle e\rangle^\bot}\max\limits_{s}\bigl[|f(x+se)-g(x+se)|^2\rho(x+se)\bigr]\, dx
+
4\int|\partial_ef - \partial_eg||f-g|\, d\mu
\\
\le
c_5(d,n)\Bigl(\int |f-g|\, d\mu\Bigr)^2,
\end{multline*}
where
Lemma \ref{lem3.1}, Theorems \ref{t1.2} and \ref{t2.1},
and estimate~(\ref{est1}) were applied in the last inequality.

Combining the above estimates, we get the bound
$$
\int\varphi(f) - \varphi(g)\, d\mu
\le c_6(d,n)\Bigl[
\|\partial_eg\|_2^{-1/(d-1)}\varepsilon^{1/(2d-2)}
+
\varepsilon^{-1/2}\|f-g\|_1
\Bigr].
$$
Taking $\varepsilon=\bigl[\|\partial_eg\|_2^{1/(d-1)}\|f-g\|_1\bigr]^{(2d-2)/d}$,
we obtain
$$
\int\varphi(f) - \varphi(g)\, d\mu
\le 2c_6(d,n)\|\partial_eg\|_2^{-1/d}\|f-g\|_1^{1/d}.
$$

In the case of an arbitrary (isotropic log-concave)
density $\rho$ on $\mathbb{R}^n$ with $n\ge2$, the estimate follows from the approximation
by the measures with densities $\rho_\varepsilon$, where
$$
\rho_\varepsilon(x):=\rho*\psi_\varepsilon\bigl((1+\varepsilon^2)^{\frac{1}{n+2}}\cdot x\bigr),
$$
$\psi$ is the density of the standard Gaussian measure on $\mathbb{R}^n$ and
$\psi_\varepsilon(x) = \varepsilon^{-n}\psi(\varepsilon^{-1}x)$.

Finally, the one-dimensional case follows from the
case $n=2$ by consideration of the product measure $\mu\otimes\mu$
and polynomials depending only on the first argument.
The lemma is proved.
\end{proof}

\begin{corollary}\label{c3.1}
Let $n, d\in\mathbb{N}$.
Then there is a constant $c(d,n)$ depending only on
$d$ and $n$ such that, for any
isotropic log-concave measure
$\mu$ on $\mathbb{R}^n$,
any pair of polynomials  $f$ and $g$
of degree $d$ on $\mathbb{R}^n$, and
any function
$\varphi\in C_0^\infty(\mathbb{R})$ with $\|\varphi\|_\infty\le1$
one has
$$
\Bigl(\int|g - \mathbb{E}g|^{1/d}d\mu\Bigr)\int\bigl(\varphi(f) - \varphi(g)\bigr)\, d\mu\le
c(d,n)\int|f-g|^{1/d}\, d\mu.
$$
\end{corollary}

\begin{proof}
We note that
$$
\int|(\nabla g, e)|^{1/d}\, d\mu
\le\biggl(\int|(\nabla g, e)|^2\, d\mu\biggr)^{1/2d}
=\|\partial_eg\|_2^{1/d}.
$$
Hence, by Lemma \ref{lem3.2}
$$
\Bigl(\int|(\nabla g, e)|^{1/d}\, d\mu\Bigr)\int\bigl(\varphi(f) - \varphi(g)\bigr)\, d\mu
\le c(d,n)\|f-g\|_1^{1/d}.
$$
Integrating in the above inequality with respect to the normalized surface
measure $\sigma_n$ on the unite sphere, we get
$$
\Bigl(\int_{S^{n-1}}\int|(\nabla g, e)|^{1/d}\, d\mu\, \sigma_n(de)\Bigr)
\int\bigl(\varphi(f)-\varphi(f)\bigr)\, d\mu
\le
c(d,n)\|f-g\|_1^{1/d}.
$$
By Fubini's theorem
\begin{multline*}
\int_{S^{n-1}}\int|(\nabla g, e)|^{1/d}\, d\mu\, \sigma_n(de)=
\int\int_{S^{n-1}}|(\nabla g, e)|^{1/d}\, \sigma_n(de)\, d\mu\\=
\int|\nabla g|^{1/d}\int_{S^{n-1}}|(e,e_1)|^{1/d}\, \sigma_n(de)\, d\mu=
c_1(d,n)\int|\nabla g|^{1/d}\, d\mu.
\end{multline*}
So, by the above equality and Theorem \ref{t1.2}, we have
$$
\Bigl(\int|\nabla g|^{1/d}\, d\mu\Bigr)\int\bigl(\varphi(f) - \varphi(g)\bigr)\, d\mu
\le c_2(d,n)\int|f-g|^{1/d}\, d\mu.
$$
Applying Theorem \ref{t1.2} again, we get
$$
\||\nabla g|\|_2^{1/d}\le C(d)\int|\nabla g|^{1/d}\, d\mu.
$$
Thus, by
Theorem \ref{t1.6} we get the desired bound.
The corollary is proved.
\end{proof}

We are now ready to prove the main result of the paper.
The key part of the proof is the application
of the localization lemma, which enables us to reduce
considerations
to a space of dimension at most $4$.

\begin{theorem}\label{t3.1}
Let $d,n\in\mathbb{N}$, $d\ge2$.
Then, there is a constant $C(d)$ depending only on~$d$ such that, for any
log-concave measure $\mu$ on $\mathbb{R}^n$ and any pair of
polynomials $f$ and $g$ of degree $d$ on $\mathbb{R}^n$,
one has
$$
\sigma_g^{1/d}\|\mu\circ f^{-1} - \mu\circ g^{-1}\|_{\rm TV}\le
C(d)\|f-g\|_2^{1/d},
$$
where
$\displaystyle \sigma_g^2:=\mathbb{D}g=\int(g-\mathbb{E}g)^2\, d\mu,\
\mathbb{E}g:=\int g\, d\mu$.
\end{theorem}

\begin{proof}
Set $R(d):=\max\limits_{n=1,2,3,4}c(d,n)$,
where $c(d,n)$ is the constant from Corollary \ref{c3.1}.
Due to Theorem \ref{t1.2} it is sufficient to prove that,
for any function $\varphi\in C_0^\infty(\mathbb{R})$ with $\|\varphi\|_\infty\le1$,
one has
\begin{equation}\label{bound1}
\Bigl(\int|g - \mathbb{E}g|^{1/d}d\mu\Bigr)\int\bigl(\varphi(f) - \varphi(g)\bigr)\, d\mu\le
R(d)\int|f-g|^{1/d}\, d\mu.
\end{equation}

First we consider the case $n\in\{1,2,3,4\}$.
Recall that for an arbitrary log-concave measure $\mu$ on $\mathbb{R}^n$
there is a nondegenerate linear mapping $T\colon\mathbb{R}^n\to\mathbb{R}^n$
such that measure $\mu\circ T^{-1}$ is isotropic.
By Corollary \ref{c3.1}, for every pair $f,g$ of polynomials of degree $d$,
we have
\begin{multline*}
\Bigl(\int |g\circ T^{-1} - \mathbb{E}(g\circ T^{-1})|^{1/d}\, d(\mu\circ T^{-1})\Bigr)
\int\bigl(\varphi(f\circ T^{-1}) - \varphi(g\circ T^{-1})\bigr)\, d(\mu\circ T^{-1})
\\
\le R(d)\int |f\circ T^{-1}-g\circ T^{-1}|^{1/d}\, d(\mu\circ T^{-1})
\end{multline*}
as functions $f\circ T^{-1}$ and $g\circ T^{-1}$ are also polynomials of degree $d$.
This implies estimate~(\ref{bound1})
for log-concave measures on $\mathbb{R}^n$ with $n\in\{1,2,3,4\}$.

Let now $n$ be an arbitrary positive integer.
Fix a convex compact set $K$,
numbers $a,b>0$, and polynomials $f, g$ of degree $d$.
Let
$$
F_1=g,\quad F_2=-g,\quad F_3=|g|^{1/d}-a,\quad F_4=b-|f-g|^{1/d}
$$
and
let $P_{F_1,\ldots, F_4}$ be the set of
all log-concave measures supported in $K$
such that
$$
\int F_i\, d\mu\ge0, \ i=1,\ldots,4.
$$
We note that the above conditions are equivalent to the following one:
$$
\int g\, d\mu=0,\quad \int |g|^{1/d}\, d\mu\ge a,\quad \int |f-g|^{1/d}\, d\mu\le b,
$$
Consider the functional $\displaystyle \Phi_{f,g}(\mu):=\int\bigl(\varphi(f) - \varphi(g)\bigr)\, d\mu$.
Note that the restriction of a polynomial to a linear subspace
will be again a polynomial (of the same degree) on this subspace.
Thus,
$\Phi_{f,g}(\mu)\le R(d)ba^{-1}$
for an arbitrary measure $\mu\in P_{F_1,\ldots, F_4}$
such that the smallest affine subspace containing the support of $\mu$
is of dimension not greater than~$4$.
By Theorem \ref{t1.7},
$$
\Phi_{f,g}(\mu)\le R(d)ba^{-1}
$$
for any measure $\mu\in P_{F_1,\ldots, F_4}$,
implying bound (\ref{bound1}) for an arbitrary log-concave
measure on $\mathbb{R}^n$ with compact support.
The general case follows by approximation.
The theorem is proved.
\end{proof}

\vskip .1in

We now briefly discuss the infinite-dimensional case.
Let $E$ be a locally convex space equipped with the Borel $\sigma$-field
and let $E^*$ be the topological
dual space to $E$.
A Radon probability measure $\mu$ on $E$ is called
log-concave (or convex) if $\mu\circ A^{-1}$ is a log-concave measure
on $\mathbb{R}^n$
for every continuous linear operator $A\colon E\to\mathbb{R}^n$.
For a Radon probability measure $\mu$ on $E$,
denote by
$\mathcal{P}^d(\mu)$ the closure in $L^2(\mu)$ of the set of all
functions of the form
$f(\ell_1,\ldots,\ell_n)$, where $n$ is an arbitrary positive integer,
$\ell_j\in E^*$ are arbitrary continuous linear functionals, and
$f$ is an arbitrary polynomial on $\mathbb{R}^n$ of degree $d$.
It is shown in \cite{ArYar} that every function from $\mathcal{P}^d(\mu)$
has a version that is a polynomial of degree $d$ in the usual algebraic sense, i.e.,
this version is of the form
$$
b_0 + b_1(x) + b_2(x, x) + \ldots + b_d(x,\ldots,x),
$$
where each $b_j(x_1,\ldots, x_j)$ is a multilinear function on $E^j$.

\begin{corollary}
Let $d,n\in\mathbb{N}$, $d\ge2$.
Then, there is a constant $C(d)$ depending only on~$d$ such that, for any
log-concave measure $\mu$ on a locally convex space $E$
and any functions $f, g\in \mathcal{P}^d(\mu)$,
one has
$$
\sigma_g^{1/d}\|\mu\circ f^{-1} - \mu\circ g^{-1}\|_{\rm TV}\le
C(d)\|f-g\|_2^{1/d},
$$
where
$\displaystyle \sigma_g^2:=\mathbb{D}g=\int(g-\mathbb{E}g)^2\, d\mu,\
\mathbb{E}g:=\int g\, d\mu$.
\end{corollary}

\vskip .1in

The author is a Young
Russian Mathematics award winner and would like to thank its sponsors and jury.

This research was supported by the RFBR
Grant 17-01-00662, the DFG through the project RO 1195/12-1 and
the CRC 1283 at Bielefeld University, and by the Foundation for the Advancement of Theoretical
Physics and Mathematics ``BASIS''.


\begin{thebibliography}{}
\bibitem{ArYar} Arutyunyan, L. M., Yaroslavtsev, I. S.:
On measurable polynomials on infinite-dimensional
spaces. Dokl. Math. 87(2), 214--217 (2013)

\bibitem{Ball} Ball, K.:
Logarithmically concave functions and sections of convex
sets in $\mathbb{R}^n$. Studia Math. 88(1), 69--84 (1988)


\bibitem{BobkIsop} Bobkov, S. G.:
Isoperimetric and analytic inequalities for log-concave probability measures.
Annals Probab. 27(4), 1903--1921 (1999)

\bibitem{BobkPoly}
Bobkov, S. G.: Remarks on the growth of $L^p$-norms of polynomials.
In: Geometric Aspects of Functional Analysis, pp.~27--35. Lecture Notes in Math. V. 1745. Springer (2000)

\bibitem{Bobk}
Bobkov, S. G.:
Some generalizations of Prokhorov's results on
Khinchin-type inequalities for polynomials.
Theory Probab. Appl. 45(4), 644--647 (2000)

\bibitem{GM}
Bogachev, V. I.: Gaussian measures. Amer. Math. Soc., Providence, Rhode Island (1998)

\bibitem{DiffMeas}
Bogachev, V. I.: Differentiable measures and the Malliavin calculus.
Amer. Math. Soc., Providence, Rhode Island (2010)


\bibitem{MeasTheor}
Bogachev, V. I.: Measure theory. V. 1. Springer, Berlin -- New York, 2007

\bibitem{BKZ} Bogachev, V. I., Kosov, E. D., Zelenov, G. I.:
Fractional smoothness of distributions of
polynomials and a fractional analog of the
Hardy--Landau--Littlewood inequality. Trans. Amer. Math. Soc.
370(6), 4401--4432 (2018)


\bibitem{Bor}
Borell, C.: Convex measures on locally convex spaces.
Ark. Math. 12, 239--252 (1974)


\bibitem{DavM}
Davydov, Y. A., Martynova, G. V.:
Limit behavior of multiple stochastic integral. Statistics and Control of Random
Processes. Nauka, Preila, Moscow. 55--57 (1987)

\bibitem{CarWr}
Carbery, A., Wright, J.: Distributional and $L^q$ norm inequalities for polynomials
over convex bodies in $R^n$.
Math. Res. Lett. 8(3), 233--248 (2001)

\bibitem{FrGue}
Fradelizi, M., Gu\'edon, O.:
The extreme points of subsets of s-concave probabilities
and a geometric localization theorem.
Discrete Comput. Geom. 31(2), 327--335 (2004)

\bibitem{KLS}
Kannan, R., Lovasz, L., Simonovits, M.:
Isoperimetric problems for
convex bodies and a localization lemma.
Discrete Comput. Geom. 13, 541--559 (1995)

\bibitem{Kos}
Kosov, E. D.:
Fractional smoothness
of images of logarithmically concave measures under polynomials.
J. Math. Anal. Appl. 462(1), 390--406 (2018)

\bibitem{Klartag} Klartag, B.:
On convex perturbations with a bounded isotropic constant.
Geometric and Functional Analysis GAFA 16(6), 1274--1290 (2006)

\bibitem{Klartag07}
Klartag, B.:
Power-law estimates for the central limit theorem for convex sets. J. Funct. Anal. 245(1), 284--310 (2007)

\bibitem{Klar-CLT}
Klartag, B.:
A central limit theorem for convex sets.
Invent. Math. 168(1), 91--131 (2007).

\bibitem{Krug} Krugova, E. P.:
On translates of convex measures. Sbornik Math. 188(2), 227--236 (1997)

\bibitem{LS}
Lovasz, L., Simonovits, M.: Random walks in a convex body and an improved
volume algorithm. Random Structures and Algorithms 4(4), 359--412 (1993)

\bibitem{NSV}
Nazarov, F., Sodin, M., Volberg, A.: The geometric Kannan--Lovasz--Simonovits
lemma, dimension-free estimates for the distribution of the values of polynomials, and
the distribution of the zeros of random analytic functions.
St. Petersburg Math. J. 14(2), 351--366 (2003)


\bibitem{NP}
Nourdin, I., Poly, G.:
Convergence in total variation on Wiener chaos.
Stochastic Processes  Appl. 123(2), 651--674 (2013)

\bibitem{Zel}
Zelenov, G. I.:
On distances between distribution of polynomials.
Theory Stoch. Processes 38(2), 79--85
(2017)

\end{thebibliography}
\end{document}